\documentclass[12pt]{article}

\usepackage[T2A]{fontenc}
\usepackage[utf8]{inputenc}

\usepackage{amsfonts}
\usepackage{amssymb}
\usepackage{amsmath}
\usepackage{amsthm}

\def \le {\leqslant}

\def \ge {\geqslant}

\theoremstyle{plain}

 scaled \magstep2

\begin{document}

\begin{Large}

\centerline{\bf Schmidt's conjecture and Badziahin-Pollington-Velani's theorem}

\vskip+0.5cm
\centerline{by    Nikolay G. Moshchevitin\footnote{research is supported by RFBF grant No. 09-01-00371a}
} 
\end{Large}
\vskip2.0cm

\begin{small}
{\bf Abstract.}\,\,\,
We give a simplified exposition of 
  the easiest case of 
a breakthrough result by D.Badziahin, A.Pollington and
S.Velani related to W.M.Schmidt's conjecture.
\end{small}
\vskip2.0cm

{\bf 1.  Schmidt's conjecture.}

In this paper all numbers are real.

For $\alpha, \beta \in [0,1] $ under the condition $ \alpha +\beta = 1 $ and $ \delta >0$
 we consider the sets
$$
{\rm BAD}(\alpha, \beta ;\delta ) = \left\{\xi = (\xi_1,\xi_2 ) \in [0,1]^2:\,\,\,\inf_{p\in \mathbb{N}} \,\, \max \{ p^\alpha ||p\xi_1||,
 p^\beta ||p\xi_2||\} \ge \delta \right\}
$$
(here $||\cdot ||$ denotes the distance to the nearest integer) and
$$
{\rm BAD}(\alpha, \beta  ) = \bigcup_{\delta > 0} {\rm BAD}(\alpha, \beta ;\delta )
 .
$$
  In \cite{SCH} Wolfgang M. Schmidt conjectured that for any $\alpha_1,\alpha_2, \beta_1,\beta_2 \in [0,1],\, \alpha_1 +\beta_1 =\alpha_2+\beta_2 =1$
  the intersection
  $$
{\rm BAD}(\alpha_1, \beta_1  ) \bigcap {\rm BAD}(\alpha_2, \beta_2  )
$$
is not empty. This conjecture was recently proved in a breakthrough paper by  Dzmitry Badziahin, Andrew Pollington and Sanju Velani \cite{PSV}.
They proved a more general result: {\it for any finite collection of  pairs $(\alpha_j,\beta_j),\,\,
0\le \alpha_j,\beta_j \le 1,\,\, \alpha_j+\beta_j = 1,\,\, 1\le j \le r$ 
and for any $\theta$ under the condition
$$
\inf_{q\in \mathbb{N}} q||q\theta|| >0
$$
the intersection
\begin{equation}\label{jjj}
\bigcap_{j=1}^r \{ \xi \in [0,1]:\,\,\,
(\theta, \xi ) \in
 {\rm BAD}(\alpha_j, \beta_j  )
\}
\end{equation}
has full Hausdorff dimension.}

Moreover one can take a certain infinite intersection in (\ref{jjj}).

This result was obtained by an original method invented by 
D.Badziahin, A.Pollington and S.Velani.
In the present paper we do not obtain  any new result.
The main purpose of the present paper is to give a more clear exposition of 
Badziahin-Pollington-Velani's method in the easiest case.

{\bf 2. The simplest case.}

 The result by  D.Badziahin, A.Pollington and S.Velani in the form (\ref{jjj}) is non-trivial  even for one set
${\rm BAD}(\alpha, \beta)$ and even for $\alpha=\beta = \frac{1}{2}$.
In this case  the result is as follows:
{\it for $\theta$  such that 
\begin{equation}\label{jjj1}
\inf_{q\in \mathbb{N}} q^2||q\theta|| >0
\end{equation}
 the set 
 $$\{ \xi \in [0,1]:\,\,\,
(\theta, \xi ) \in
 {\rm BAD}(1/2,1/2  )
\}
$$
has full Hausdorff dimension.}

In the dual form the result proclaims that under the condition (\ref{jjj1}) the set
$$
\{\xi \in [0,1]:\,\,  \inf_{(A,B) \in \mathbb{Z}^2 \setminus \{(0,0)\}} ||A\theta -B\xi ||\cdot \max(A^2,B^2) >0\}
$$
 has full Hausdorff dimension.

In the present paper we show how Badziahin-Pollington-Velani's 
construction 
gives a proof of the following result.

{\bf Proposition 1.}\,\,{\it 
Let 
\begin{equation}\label{000}
0<\delta \le 2^{-1622}
.
\end{equation}
 Suppose that
\begin{equation}\label{0}
\inf_{q\in \mathbb{N}} q^2||q\theta|| \ge\delta.
\end{equation}
Then there exists $\xi$ such that for all integers $A,B$ with $\max (|A|,|B|) >0$  one has
\begin{equation}\label{00}
 ||A\theta -B\xi ||\cdot \max(A^2,B^2) \ge {\delta} .
\end{equation}
}

Of course the constant $2^{1622}$ in   (\ref{000}) may be reduced.

In sections 4 -  10 we give a complete proof of Proposition 1.

{\bf 3. Acknowledgements.}

  First of all the author thanks Dzmitry Badziahin for visiting Moscow in January 2010. His explanations were 
very useful for the author for undestanding the  construction of the proof.  Then the author would like to thank the
University of Astrakhan as this paper was written mostly during the author's visit of the University of Astrakhan.
The author thanks all the participants of the Astrakhan's seminar in Number Theory and especially Renat Akhunzhanov.
The author thanks all the participants of number-theoretical seminars in the Lomonosov Moscow State University and in the
Independent University of Moscow for useful discussions. Especially the author would like to thank Irog Rochev
for his suggestions related to 
an earlier version of the present manuscript.

{\bf 4. Parameters.}

Suppose
$$R\ge 2^{422}$$
to be an integer.  The integer parameter  $ n$ increases to $+\infty$.
Let
$$
0<
\delta <
\frac{1}{3R^{\frac{2533}{660}}}.
$$
Put
\begin{equation}\label{lala}
\lambda  = \frac{1741}{330},
\end{equation}

\begin{equation}\label{dk}
\kappa =\delta R^{\frac{6}{5}},
\end{equation}
so
\begin{equation}\label{00k}
\kappa \le \frac{ 1}{3R^{\frac{\lambda}{2}}}.
\end{equation}
Let
$k$ 
be an integer under the condition
\begin{equation}\label{kkk}
1\le 2^k \le R.
\end{equation}
So
\begin{equation}\label{log}
0\le k \le \left\lfloor
\frac{\log R}{\log 2}\right\rfloor.
\end{equation}
Given $k$  we define
\begin{equation}\label{de}
d_k = 
\left\lfloor
\left(
\frac{\kappa}{\delta}
\cdot \frac{2^k}{R}\right)^{\frac{2}{3}}
\cdot R^{\frac{2}{165}}
\right\rfloor
\end{equation}
and
\begin{equation}\label{ke}
K_k = 
\frac{\delta}{\kappa}
\cdot \frac{R^2}{2^k}.
\end{equation}
One  can easily see that
\begin{equation}\label{de1}
2\le \lfloor R^{\frac{8}{55}}\rfloor\le
d_k 
\le R^{\frac{134}{165}}
 \end{equation}
 and
\begin{equation}\label{kd}
d_k\cdot K_k \le
  R^{\frac{52}{55}}.
\end{equation}

{\bf 5. Lines and forbidden intervals.}

Given integers $A,B,C$ with $ (A,B,C) =1, B>0$
defne $L = L(A,B,C)$   to be a line
$$
L=L(A,B,C) = \{ (x,y)\in \mathbb{R}^2:\,\, Ax-By+C=0\}.
$$
Put
\begin{equation}\label{ah}
H(A,B) = B\max (A^2,B^2).
\end{equation}

Let 
$$
\Delta = \Delta (A,B,C ) = \left( \frac{A\theta +C}{B} -\frac{\delta}{H(A,B)},
\frac{A\theta +C}{B} +\frac{\delta}{H(A,B)}\right)
$$
be the interval of the length
$$
|\Delta (A,B,C ) | = \frac{2\delta}{H(A,B)}.
$$
(Everywhere in the sequel $|J|$ stands for the length of an interval or a segment $J$.)

For our purpose it is enough to prove that
$$
[0,1]\setminus \left( \bigcup_{A,B,C} \Delta (A,B,C ) \right) \neq \varnothing 
$$
where the
union is taken over all
triples of integers $A,B,C$
such that
$$
B>0,\,\,\, (A,B,C) = 1.
$$

It is convenient to consider the segment
$$
\Theta = \{ (x,y)\in \mathbb{R}^2:\,\, x = \theta , \,\, 0\le y \le 1\}
.
$$
Also it is convenient to consider intervals
$$
\overline{\Delta} = \overline{\Delta} (A,B,C )
=\{ (x,y)\in \mathbb{R}^2:\,\, x = \theta , \,\,  y\in \Delta (A,B,C)\}
 .
$$
So our task is to prove that
\begin{equation}\label{iii}
\Theta \setminus \left( \bigcup_{A,B,C} \overline{\Delta} (A,B,C ) \right) \neq \varnothing .
\end{equation}

{\bf 6. Inductive construction.}

We describe an inductive procedure  to establish (\ref{iii}).  We take an arbitrary segment
$$J_1 \subset \Theta
$$
of the length $$ |J_1|=\frac{\kappa}{R}.$$ Now we describe the inductive process of constructing segments $J_n^\nu, n=1,2,3,...$.

Given an integer  $n\ge 1$ suppose we have a non-empty collection of segments 
\begin{equation}\label{coll}
J_n^\nu\subset \Theta,\,\,\,\, 1\le \nu \le T_n,
\end{equation}
$$
|J_n^\nu|=\frac{\kappa}{R^n}
$$
such that
\begin{equation}\label{inductive}
J_n^\nu \cap \overline{\Delta }(A,B,C) =\varnothing
\end{equation}
for all triples $A,B,C$ under consideration such that $ H(A,B) < R^{n-1}$.
(For $n=1$ this condition is empty.)

Each of the segments $J_n^\nu$ we divide into $R$ equal segments
\begin{equation}\label{form}
I_{n+1}^{\nu,\mu},\,\,\,\, 1\le \mu \le R,
\end{equation}
so
$$
J_n^\nu = \bigcup_{1\le \mu \le R} I_{n+1}^{\nu,\mu},\,\,\,\,
|I_{n+1}^{\nu,\mu } | = \frac{|J_n^\nu |}{R} = \frac{\kappa}{R^{n+1}}.
$$ 

We must consider the collection of the intervals 
\begin{equation}\label{step}
\overline{\Delta} = \overline{\Delta} (A,B,C ),
 \end{equation}
\begin{equation}\label{step1}
R^{n-1} \le H(A,B)< R^n
\end{equation}
and prove that   among the segments
\begin{equation}\label{segments}
I_{n+1}^{\nu,\mu},\,\,\,\, 1\le \mu \le R,\,\,\, 1\le \nu \le T_n
\end{equation}
there exist a large number of segments $
I_{n+1}^{\nu,\mu}$ such that
\begin{equation}\label{cappo}
I_{n+1}^{\nu,\mu}\cap \overline{\Delta}(A,B,C) = \varnothing
\end{equation}
for all intervals $\overline{\Delta}$   of the form (\ref{step}) satisfying (\ref{step1}) .

In order to do this 
for any natural $m\le n$ we must consider the corresponding collection
of the segments
$$
J_m^\nu\subset \Theta,\,\,\,\, 1\le \nu \le T_m,\,\,\,\,\,
|J_m^\nu|=\frac{\kappa}{R^m}
$$
such that
$$
J_m^\nu \cap \overline{\Delta }(A,B,C) =\varnothing
$$
for all triples $A,B,C$ under consideration such that $ H(A,B) < R^{m-1}$.

Obviuosly $T_n>0$ implies $T_m > 0$ for all $m\le n$ as the collections are nested:
$$
\bigcup_{1\le \nu \le T_1}J_1^\nu\supset \cdots \supset 
\bigcup_{1\le \nu \le T_m}J_m^\nu\supset
\bigcup_{1\le \nu \le T_{m+1}}J_{m+1}^\nu\cdots\supset
\bigcup_{1\le \nu \le T_n}J_n^\nu
.
$$

It happens that in order to show that many segments of the form (\ref{segments})
 satisfy (\ref{cappo}) we
must
 assume that for any natural $m \le n$ we have a certain lower bound for the quantity
$T_m$.
All precise estimates and inequalities will be formulated in the next sections.

{\bf 7. Single interval $\overline{\Delta}(A,B,C)$.}

Remind that the interval $\overline{\Delta}(A,B,C)$ has the length equal to $
|\overline{\Delta} (A,B,C ) | = \frac{2\delta}{H(A,B)} 
$. So given   $\overline{\Delta}(A,B,C) $ the number of  segments $
I_{n+1}^{\nu,\mu}$ satistying 
\begin{equation}\label{uuuu}
I_{n+1}^{\nu,\mu}\cap \overline{\Delta}(A,B,C) \not= \varnothing
\end{equation}
is 
\begin{equation}\label{bou}
\le \frac{|\overline{\Delta}(A,B,C)|}{ |I_{n+1}^{\nu,\mu}|}+2 =
\frac{2\delta}{\kappa}\cdot \frac{R^{n+1}}{H(A,B)}+2.
\end{equation}

Given $k$ from the interval (\ref{kkk}) consider the following condition on $H(A,B)$ which is  stronger than the condition   (\ref{step1}):
\begin{equation}\label{ennk}
2^kR^{n-1}\le H(A,B)  = B\cdot \max (A^2,B^2) < \min ( 2^{k+1}R^{n-1};R^n).
\end{equation}

Let $A,B$ satisfy the condition (\ref{ennk}).
Consider a fixed interval $\overline{\Delta}(A,B,C).$
We see 
(here we should refer to the definition (\ref{ke}) of the parameter $K_k$)
 that the the number of segments $
I_{n+1}^{\nu,\mu}$ satisfying (\ref{uuuu}) with fixed $A,B,C$ is less or equal than
\begin{equation}\label{val}
2K_k+2.
\end{equation}.

{\bf 8.  Lines with bounded coefficient $|A|/B$.}

In this section we consider a single segment $J_n = J_n^\nu$ from the collection (\ref{coll}). 

Given $k$ from the interval (\ref{kkk}) consider  all the lines $L (A,B,C)$ such that  coefficients $A,B$ satisfy the condition
(\ref{ennk})
 and the {\it additional } condition
\begin{equation}\label{ennadd}
B > R^{\frac{n}{3}-\lambda}.
\end{equation}
 
The purpose of
  the current section is to prove that 
{\it   the number of   segments of the form
(\ref{form}) satisfying
 (\ref{uuuu})
for some interval  $\overline{\Delta} (A,B,C)$ under conditions (\ref{ennk},\ref{ennadd})
is  $$\le \gamma R^{\frac{52}{55}} .$$}

An admissible value for $\gamma$ is $\gamma = 2^{13}$.
 
Recall that  $k$
satisfies (\ref{log}). So from the desired upper bound 
for the number of segments satisfying (\ref{uuuu})
we see  that {\it   the number of   segments of the form
(\ref{form}) satisfying
 (\ref{uuuu})
for some  interval  $\overline{\Delta} (A,B,C)$ under conditions (\ref{step1},\ref{ennadd})
is  
$$
\le \gamma R^{\frac{52}{55}} \left(\left\lfloor\frac{\log R}{\log 2}\right\rfloor+1\right).$$}

Note that under the conditions (\ref{ennk},\ref{ennadd})
one has
\begin{equation}\label{B}
B\le 2^{\frac{k+1}{3}}R^{\frac{n-1}{3}} 
\end{equation}
and
\begin{equation}\label{A}
|A|\le 2^{\frac{1}{6}}\cdot 2^{\frac{k+1}{3}}R^{\frac{n-1}{3}+\frac{\lambda}{2}} .
\end{equation}
The last inequality follows from 
$$
A^2 \le
\frac{
2^{k+1}R^{n-1}}{B} <
2^{k+1}R^{\frac{2(n-1)}{3}-\frac{1}{3}+\lambda} \le
2^{\frac{1}{3}}\cdot 2^{\frac{2}{3}(k+1)}
R^{\frac{2}{3}(n-1)+\lambda}.
$$

Also from (\ref{ennk},\ref{ennadd} ) we see that

$$
\left(\frac{|A|}{B}\right)^2\le \frac{H(A,B)}{B^3} <\frac{R^n}{R^{n-3\lambda}} = R^{3\lambda}.
$$
So
\begin{equation}\label{AB}
\frac{|A|}{B} \le  R^{\frac{3\lambda}{2}}. 
\end{equation}

{\bf 8.1. Lemmata about lines intersecting a segment.}

Here we give few lemmas. They will be useful not only in Section 8  but also in Section 9
where we consider a general situation.

{\bf Lemma 1.}\,\, {\it Consider a segment  $J_n = J_n^\nu$ from the collection
(\ref{coll}). Suppose that there exist two lines
$$
L_1= L(A_1,B_1,C_1),\,\,\, L_2= L(A_2,B_2,C_2)
$$
such that
$$
L_i \cap J_n \neq \varnothing ,\,\,\,\, i = 1,2
$$
and
$$
H(A_1,B_1), H(A_2,B_2) < R^n.
$$
 Then lines $L_1$ and $L_2$ are not parallel.}

Proof.  \,\, Lines $L_i, i = 1,2$  intersect the segment $J_n \subset \Theta$ in 
 points 
$$
\left( \theta, 
\frac{A_1\theta+C_1}{B_1}\right),\,\,\,\
\left(\theta,
\frac{A_2\theta+C_2}{B_2}\right)
$$
with $y$-coordinates
$$
\frac{A_1\theta+C_1}{B_1},\,\,\,\
\frac{A_2\theta+C_2}{B_2},\,\,\,\,
\left|\frac{A_1\theta+C_1}{B_1}-
\frac{A_2\theta+C_2}{B_2}\right|\le |J_n|.
$$
Suppose these lines to be parallel.
Then
$$
\frac{A_1}{B_1} =\frac{A_2}{B_2}
$$
and
$$
\frac{1}{B_1B_2} \le 
\left|
\frac{C_1}{B_1} -\frac{C_2}{B_2}
\right| =
\left|
\frac{A_1\theta+C_1}{B_1} -
\frac{A_2\theta+C_2}{B_2}\right| \le |J_n| = \frac{\kappa}{R^n}.
$$
From the inequality $ B^3_j \le H(A_j,B_j) < R^n$ 
we see that $B_j < R^{\frac{n}{3}}$ and so
$$
\frac{1}{B_1B_2}\ge \frac{1}{R^{\frac{2}{3}n}}.$$
 As $\kappa$ is small enough we have a contradiction. $\Box$

{\bf Lemma 2.}\,\, {\it Consider a segment  $I\subset \{(x,y):\,\, x =\theta\}$ of the length $|I|$.
Suppose that two lines
$$
L_1= (A_1,B_1,C_1),\,\,\, L_2= L(A_2,B_2,C_2)
$$
intersect this segment $I$.
Suppose that
$$
L_1\cap L_2 = P = \left(\frac{p}{q},\frac{r}{q}\right),\,\,\,
p,r,q\in \mathbb{Z},\,\,\, q>0, \,\,\, (p,r,q) = 1.
$$
Then

{\rm (i)}\,  $ |q\theta - p| \le |I| B_1B_2$;

{\rm (ii)}\,  $ q \le 2 \,\max (|A_1|, |A_2|) \, \max (B_1, B_2)$.}

Proof.
\,\,
Obviously  rational numbers $\frac{p}{q},\frac{r}{q}$
 satisfy
$$
A_i \frac{p}{q}- B_i \frac{r}{q}+C_i = 0,\,\,\, i = 1,2.
$$
So
$$
\frac{p}{q}
=\frac{B_1C_2 - B_2C_1}{A_1B_2-A_2B_1},\,\,\,
\frac{r}{q} = \frac{A_1C_2 - A_2C_1}{A_1B_2-A_2B_1}.
$$
As $(p,r,q) = 1$ we see that for some non-zero integer $s$ one has
$$
sq = A_1B_2 - A_2B_1,
$$
$$
sp = B_1C_2 - B_2C_1,
$$
$$
sr = A_1C_2 - A_2C_1.
$$
So (ii) follows from the first of these three equalities as
$$
q \le |s|q  = |A_1B_2 - A_2B_1| \le  2\max (|A_1|, |A_2|) \, \max (B_1, B_2).
$$
Now
$$
|q\theta - p| \le |sq\theta - sp| =
B_1B_2 \cdot \left | \frac{A_1\theta+C_1}{B_1} - 
\frac{A_2\theta+C_2}{B_2}\right| \le
B_1B_2 |I|,
$$ and (i) follows. $\Box$ 

{\bf Lemma 3.}\,\, {\it
All the lines $ L = L(A,B,C)$ 
such that 
$$
L(A,B,C) \cap J_n \neq \varnothing,
$$
$$
H(A,B) < R^n
$$
satisfying the additional condition (\ref{ennadd})
have a single common point.}

Proof.\footnote{This  proof was  suggested to the author by Igor Rochev.}

\,\, From Lemma 1 it follows that any two lines  intersecting $J_n$  have a common point.
Suppose that we have three lines 
$$
L_i=
L(A_i,B_i,C_i),\,\,\,\, i = 1,2,3
$$
intersecting $J_n$  which satisfy the conditions of Lemma 3 but do not have  a common point.
Then
$$
D =
\left|
\begin{array}{ccc}
A_1&B_1&C_1\cr
A_2&B_2&C_2\cr
A_3&B_3&C_3
\end{array}
\right|\neq 0.
$$
Let 
$(\theta, \xi) $ be the middle point of the segment $J_n$.
Then
$$
1\le |D| =|\,
\left|
\begin{array}{ccc}
A_1&B_1&A_1\theta -B_1\xi +C_1\cr
A_2&B_2&A_2\theta -B_2\xi +C_2\cr
A_3&B_3&A_3\theta -B_3\xi +C_3
\end{array}
\right|\, |.
$$
 Suppose that $(\theta, Y_i) = J_n\cap L_i$. Then
$$
Y_i = \frac{A_i\theta+C_i}{B_i}.
$$
Now
$$
|A_i\theta -B_i\xi +C_i|=
B_i |\xi - Y_i| \le  B_i \cdot \frac{\kappa}{2R^n}.
$$
Define $$A = \max_{i=1,2,3} |A_i|,\,\,\, B = \max_{i=1,2,3} B_i.$$
Then
$$
1\le |D| \le 3\cdot 2AB \cdot \frac{\kappa B}{2R^n} = \frac{3\kappa AB^2}{R^n}.
$$
We have $ B< R^{\frac{n}{3}}$
and
$|A_i|^2 < \frac{R^n}{B_i} < R^{\frac{2}{3} n +\lambda}.$  Recall that we suppose the condition  (\ref{00k}) to be valid. So
$$
1\le |D| < 3\kappa R^{\frac{\lambda}{2}} \le 1.$$ 
This is not possible and lemma is proved.$\Box$

{\bf Lemma 4.}\,\,\,{\it  Consider two lines $L_i = L(A_i,B_i,C_i), i = 1,2$.
Suppose that $A_1,B_1$ satisfy the condition (\ref{ennk}).
Suppose that 
$$
B_1\ge B_2.
$$
Let 
$$
L_1\cap L_2 = P=\left(\frac{p}{q},\frac{r}{q}\right).
$$
Put
\begin{equation}\label{pointsp}
  Y_i = \frac{A_i\theta +C_i}{B_i}
\end{equation}
and suppose that
\begin{equation}\label{igrek}
|Y_1 - Y_2| \le \frac{|J_n|}{d_k} = \frac{\kappa}{d_kR^n}.
\end{equation}
Put
\begin{equation}\label{ceka}
\sigma_k = \frac{2\kappa}{d_k}\cdot \frac{2^k}{R} \cdot \frac{1}{|q\theta-p|}. 
\end{equation}
Then
\begin{equation}\label{ppo}
B_1\le \sigma_k,\,\,\, A_1^2\le \sigma_k B_1.
\end{equation}}

Proof.\,\,\,
From (\ref{igrek}) it follows that lines $L_1,L_2$ intersect 
the line $ \{(x,y):\,\, x = \theta\}$
 in points
\begin{equation}\label{points}
{\cal Y}_1= (\theta, Y_1,),\,\,\,{\cal Y}_2=(\theta ,Y_2)
\end{equation}
where $Y_i$ are defined in (\ref{pointsp}).
We apply statement (i) of Lemma 2 with respect to the segment  $I =[{\cal Y}_1,{\cal Y}_2]$ of the length
$|I| \le \frac{\kappa}{d_kR^n}$
 to obtain the inequality
$$
|q\theta - p| \le \frac{\kappa}{d_kR^n}\cdot B_1B_2
\le
\frac{\kappa}{d_kR^n}\cdot B_1^2.
$$
From the condition (\ref{ennk}) we see that
$$
H(A_1,B_1) \le 2^{k+1}R^{n-1},
$$
and hence
$$
R^n \ge \frac{1}{2} \cdot \frac{R}{2^k} \cdot H (A_1,B_1).
$$
So
$$
|q\theta  -p| \le
\frac{2\kappa}{d_k}\cdot \frac{B_1^2}{H(A_1,B_1)}\cdot \frac{2^k}{R}
$$
or
$$
\max\left(
\frac{A_1^2}{B_1}, B_1
\right)=
\frac{H(A_1,B_1)}{B_1^2} \le
\frac{2\kappa}{d_k}\cdot \frac{2^k}{R}\cdot\frac{1}{|q\theta -p|} =\sigma_k
$$
 and Lemma 4 follows.$\Box$

{\bf  8.2. Technical lemma.}

In this section we prove a statement concerning the maximal value of the quantity $|A|/B$ under certain conditions.

{\bf Lemma 5.}\,\,\,{\it  
Let $\sigma ,W>0$.
Suppose that real numbers $A,B$   
satisfy the following conditions:
$$
0<B\le \sigma,\,\,\,\, A^2\le \sigma B,\,\,\,\, H(A,B) \ge W.
$$
Then
$$
\frac{|A|}{B} \le \left(\frac{\sigma^3}{W}\right)^{\frac{1}{4}}
.$$}

Proof.  \,\,\, Obviously the maximal value of the ratio $|A|/B$ occurs at the point $(A_*,B_*)$ which is a solution of the system
$$
\begin{cases}
A^2= \sigma B, \cr
A^2B = W.
\end{cases}
$$
So
$$
|A_*| =(\sigma W)^{\frac{1}{4}},\,\,\,\, B_* =\left(\frac{W}{\sigma}\right)^{\frac{1}{2}},
$$ 
and Lemma 5 follows. $\Box$

{\bf Collections  $\hbox{\got A}$   and $\hbox{\got B} $.}
   
Let
\begin{equation}\label{lines}
L_1,L_2,..., L_M
\end{equation}
be all the lines  $L(A,B,C)$
under conditions (\ref{ennk},\ref{ennadd}) intersecting the segment $J_n$.
Suppose that $ M \ge 2$.
From  Lemma 3 we know that all these lines  pass through a single rational point
$$P=\left(\frac{p}{q},\frac{r}{q}\right)
=\bigcap_{1\le i\le M}L_i.
$$
Put
$$
W_k= 2^{k}R^{n-1}, \,\,\, V_k = \left(\frac{\sigma_k^3}{W_k}\right)^{\frac{1}{4}}$$
(here $\sigma_k$ is defined in (\ref{ceke})) and
\begin{equation}\label{amega}
\omega_k
 = \left|\theta -\frac{p}{q}\right| \cdot  V_k  
 = \left|\theta -\frac{p}{q}\right| \cdot
\left(\frac{\sigma_k^3}{R^n} \frac{R}{2^k}\right)^{\frac{1}{4}}=
\frac{(2\kappa )^{\frac{3}{4}}}{d_k^{\frac{3}{4}} qR^{\frac{n}{4}}} |q\theta - p|^{\frac{1}{4}}
\left(\frac{2^k}{R}\right)^{\frac{1}{2}}
.
\end{equation}

We divide the collection of all the lines (\ref{lines})  into two subcollections 
$\hbox{\got A}$   and $\hbox{\got B} $.

Suppose that the collection $\hbox{\got A}$  consist of all lines of the form 
(\ref{lines}) that intersect  the segment
\begin{equation}\label{ilines}
\Omega = \left\{ (x,y) \in \mathbb{R}^2: \,\,\, x= \theta ,\,\,\,
y\in 
\left[ \frac{r}{q} - \omega_k, \frac{r}{q}+\omega_k \right] \right\}\subset
\{ (x,y):\,\, x= \theta\}.
\end{equation}
Suppose that the collection 
$\hbox{\got B} $ consists of  
all lines of the form 
(\ref{lines}) that do not intersect  the interval (\ref{ilines}).

{\bf Lemma 6.}\,\,\,{\it The number of elements in the collection $\hbox{\got B} $ is bounded by
$$
\#\hbox{\got B} 
\le d_k
.
$$}

Proof. \,\,\,Suppose that 
$\#\hbox{\got B} 
> d_k
.
$
Then there exist two lines $L_1= L(A_1,B_1,C_1),L_2 = L(A_2,B_2,C_2) \in \hbox{\got B} 
$ such that for the points
(\ref{points}) the inequality (\ref{igrek}) is valid.
Without loss of generality assume that $B_1\ge B_2$.
So we can apply Lemma 4 to see that $A_1,B_1$ satisfy  inequalities (\ref{ppo}).
It means that $A_1,B_1$  satisfy the conditions of Lemma 5 with $\sigma = \sigma_k, W = W_k$.
From Lemma 5 it follows that 
$$
\frac{|A_1|}{B_1} \le V_k.
$$
From the definition of $\omega_k$ we see  that
$$
Y_1 = \frac{A_1\theta+C_1}{B_1} \in \left[ \frac{r}{q} - \omega_k, \frac{r}{q}+\omega_k \right].
$$
So
$$
L_1 \cap \Omega = \left( \theta, Y_1\right) \neq \varnothing.
$$
It means that $L_1\in \hbox{\got A }$. This is a contradiction. $\Box$

In next two sections we deal with the collection $\hbox{\got A }$.

{\bf 8.3. Collection  $\hbox{\got A }$: the first principal inequality.}

We suppose that
\begin{equation}\label{agret}
\# \hbox{\got A } \ge 2d_k.
\end{equation}
Under this condition we deduce {\bf the first principal inequality}:

{\bf Lemma 7.}\,\,\,{\it
Suppose that (\ref{agret}) is valid. Then
$$
qd_k \le {12\sigma_k^2} .
$$}

Proof. \,\,\,  We divide the interval $J_n$ into $d_k$   intervals  $J_n(\mu)$ of the equal length
$$
|J_n(\mu )| =\frac{|J_n |}{d_k} =\frac{\kappa}{d_kR^n}.
$$
Given interval $J_n(\mu )$ consider a single  line $L=L(A,B,C)= L(\mu) $ from the collection (\ref{lines})
intersecting $J_n(\mu )$ and such that the coefficient $B$ is the smallest one.
(Of course such a line exists only in the case when the set of lines of the form  (\ref{lines}) intersecting $J_n(\mu )$
is not empty.)  So the number of lines
$L(\mu)$ is bounded by the number of intervals $J_n(\mu )$, that is $d_k$.
From Lemma 4 we see that for  any line $L(A,B,C)$ from the collection 
$\hbox{\got A } $ different from lines $L(\mu)$  its coefficients must satisfy 
(\ref{ppo}) and hence
\begin{equation}\label{set}
\max (|A|, B) \le \sigma_k.
\end{equation}
We see that there exist  $\ge d_k$ different lines from the  collection $\hbox{\got A }$ 
with coefficients satisfying (\ref{set}).
Now we should make two observations.

1. All the lines from the collection $\hbox{\got A } $ pass through the rational point
$P =\left(\frac{p}{q},\frac{r}{q}\right)$.
So the  corresponding integer points $(A,B)$
must belong to  the lattice
$$
\Lambda 
= \{ (A,B) \in \mathbb{Z}^2:\,\,\, Ap-Br \equiv 0\pmod{q}\}
$$
with the fundamental determinant
$$
{\rm det }\Lambda = q.
$$

2. As there is no parallel lines in the collection $\hbox{\got A } $
(Lemma 1)
we see that the convex hull
$$
\Pi = {\rm conv} \left( \{(0,0)\} \cup  
\{(A,B) :\,\, \exists\, C\,\,\, L(A,B,C ) \in \hbox{\got A }, \,\, L(A,B,C)\neq L(\mu )\}\right)
 $$
is a polygon with positive measure ${\rm mes} \, \Pi$ (the last inequality takes into account that $d_k \ge 2$).
We see that $\Pi$ contains $> \#\hbox{\got A }  -d_k\ge d_k$  points of the lattice $\Lambda$ (here we make use of the condition (\ref{agret})).

As the fundamental determinant of the lattice $\Lambda$ is equal  to $q$, by Pick's formula
we have 
\begin{equation}\label{ooo}
{\rm mes} \,\Pi>
\frac{q(\#\hbox{\got A }  -d_k)}{6}\ge 
\frac{qd_k}{6}.
\end{equation}
But from (\ref{set}) it follows that
$$
\Pi \subset \{(A,B)\in \mathbb{R}^2:\,\,\, \max (|A|,B) \le \sigma_k,\,\, B\ge 0\}.
$$
So
\begin{equation}\label{ooo0}
{\rm mes}\, \Pi
\le 2\sigma_k^2.
 \end{equation}
Lemma 7 immediately follows from (\ref{ooo},\ref{ooo0}).$\Box$

{\bf Lemma 8.}\,\,\,{\it  Under conditions of Lemma 7  one has
$$
|q\theta-p|\le \frac{2\sqrt{12}\,\kappa}{d_k^{\frac{3}{2}}q^{\frac{1}{2}}}\cdot \frac{2^k}{R}.
$$
}

Proof.\,\,\, Lemma 8 follows immediately from Lemma 7 and the definition of $\sigma_k$ (equality (\ref{ceka})).
$\Box$

{\bf 8.4. Interval $\Omega$: the second principal inequality.}

If $ \hbox{\got A }\neq \varnothing$ then
\begin{equation}\label{NON}
\Omega \cap J_n \neq \varnothing.
\end{equation}
This fact leads to the {\bf second principal inequality}:

{\bf Lemma 9.}\,\,\,{\it
Suppose that  under the conditions of Lemma 7 one has
\begin{equation}\label{qu}
q< R^{\frac{2}{3} (n-1)}.
\end{equation}
Then
 for the value $\omega_k$ defined in (\ref{amega}) one has
$$
\omega_k \ge \frac{\delta}{2q^{\frac{3}{2}}}.
$$}

Proof.

We apply pigeonhole principle to see that there exist integers $A,B,C$ such that $ (A,B)\neq (0,0)$ and
$$
Ap-Br+Cq=0
$$
and
\begin{equation}\label{maa}
\max (|A|,B) \le q^{\frac{1}{2}},\,\,\, B\ge 0.
\end{equation}

In fact we prove that $B>0$. Indeed if $B=0$ then $A\neq 0$ and 
one has
$$
|A\theta+C| =\left|A\left(\theta -\frac{p}{q}\right)  +A\frac{p}{q}+C\right|=
\left|A\left(\theta -\frac{p}{q}\right)\right|
= \frac{|A|}{q} |q\theta - p| \le
\frac{|q\theta-p|}{q^{\frac{1}{2}}}.
$$
From Lemma 8  and  (\ref{de}) we see that
$$
|A\theta+C|\le 
 \frac{2\sqrt{12}\,\kappa}{d_k^{\frac{3}{2}}q }\cdot \frac{2^k}{R}
<
\frac{4\sqrt{12}\,\delta}{ q  R^{\frac{1}{55}}}\le \frac{\delta}{q}
  $$
(as $R^{\frac{1}{55}}\ge 4\sqrt{12}$).
But from (\ref{0}) we see that
$$
|A\theta +C | \ge \frac{\delta}{A^2} \ge \frac{\delta}{q}.
$$
So we have a contradiction and hence $B>0$.

Now
$$
\left| A\theta -B\frac{r}{q} +C\right|=
\left| A\left(\theta-\frac{p}{q}\right) +\frac{Ap-Br+Cq}{q}\right| =
\left| A\left(\theta-\frac{p}{q}\right) \right|\le \frac{|q\theta- p|}{q^{\frac{1}{2}}},
$$
or
\begin{equation}\label{YYY}
 \left| 
\frac{A\theta+C}{B} -\frac{r}{q}\right|\le\frac{|q\theta- p|}{Bq^{\frac{1}{2}}}\le
 \frac{2\sqrt{12}\,\kappa}{B\cdot d_k^{\frac{3}{2}}q  }\cdot \frac{2^k}{R}.
\end{equation}
(here we apply Lemma 8 again).

The number $\frac{A\theta+C}{B} $ corresponds to the center of the interval
$$
\overline{\Delta}(A,B,C)
$$
with
$$
H(A,B) \le  (\sqrt{q})^3  <R^{n-1} 
$$
(here we make use of (\ref{qu},\ref{maa})).
From our  inductive assumption   in this situation one has
$$
J_n \cap  \overline{\Delta}(A,B,C) =\varnothing
$$
(see (\ref{inductive})).
Let
$$
{\cal Y }=\left(\theta , \frac{A\theta+C}{B}\right)
$$
be the center of the interval $ \overline{\Delta}(A,B,C)$. One has
\begin{equation}\label{YYYY}
{\rm dist}  (J_n , {\cal Y} ) \ge \frac{|\overline{\Delta}(A,B.C)|}{2} = \frac{\delta}{H(A,B)} =
\frac{\delta}{B \cdot \max (A^2,B^2)}\ge \frac{\delta}{Bq}.
\end{equation}
Note that the point
$$
{\cal Y}_* = \left(\theta, \frac{r}{q}\right)
$$
is the center of the segment $\Omega$.
From (\ref{NON}) it follows that
\begin{equation}\label{YYYYY}
{\rm dist}  (J_n , {\cal Y}_* )
\le \omega_k.
 \end{equation}

Now we collect together
(\ref{YYY},\ref{YYYY},\ref{YYYYY}) to see that
$$
\omega_k + 
  \frac{2\sqrt{12}\,\kappa}{B\cdot d_k^{\frac{3}{2}}q  }\cdot \frac{2^k}{R} \ge
\frac{\delta}{Bq}.
$$
But as  $R^{\frac{1}{55}}\ge 8\sqrt{12}$
we see that
$$
\frac{2\sqrt{12}\,\kappa}{  d_k^{\frac{3}{2}} }\cdot \frac{2^k}{R} \le
\frac{\delta}{2}.
$$
So
$$
\omega_k   \ge
\frac{\delta}{2Bq}\ge \frac{\delta}{2q^{\frac{3}{2}}}.
$$
Lemma 9 is proved.$\Box$

{\bf 8.5. The first fundamental lemma.}

{\bf Fundamental Lemma 1.}\,\,\,{\it
Suppose we have a segment $J_n = J_n^\nu$ satisfying (\ref{inductive}).
Then the number of segments $I_{n+1}^{\nu,\mu}$ of the form (\ref{form}) which  has non-empty intersection  with some  interval
$$
\overline{\Delta}(A,B,C)
$$
with $A,B$ satisfying (\ref{step1}) and (\ref{ennadd}) is
$$
\le 2^{13} R^{\frac{52}{55}}\log R.
$$
}

Proof.

1. Consider all values of parameter $k$ for which $ M < 3 d_k$.
For these $k$ one can see that the number of lines from  (\ref{lines})
intersecting $J_n$ is less than $3d_k$. For each line from (\ref{lines}) the corresponding interval
$\overline{\Delta}(A,B,C)$ can intersect not more than $2K_k+2$ segments of the form  (\ref{form}).
It may happen that a line $L(A,B,C)$ does not intersect the segment $J_n$ but the corresponding interval
$\overline{\Delta}(A,B,C)$  does intersect. But obviously such intervals can totally intersect not more than
$2K_k+2$ segments of the form  (\ref{form}).
So  for the parameter $k$ under consideration the number of intersected segments of the form  (\ref{form})
is 
$$
\le (2K_k+2)\cdot  3d_k\le 8R^{\frac{52}{55}}
$$
(we take into account (\ref{de1},\ref{kd})).
 
2. Consider all values of parameter $k$ for which $ M \ge  3 d_k$.
In this case we have (\ref{agret}). So Lemma 8 gives the inequality
\begin{equation}
|q\theta-p|\le \frac{2\sqrt{12}\,\kappa}{d_k^{\frac{3}{2}}q^{\frac{1}{2}}}\cdot \frac{2^k}{R}.
\label{q1}
\end{equation}
Recall that (\ref{amega}) gives
 $$
\omega_k
   =
\frac{(2\kappa )^{\frac{3}{4}}}{d_k^{\frac{3}{4}} qR^{\frac{n}{4}}} |q\theta - p|^{\frac{1}{4}}
\left(\frac{2^k}{R}\right)^{\frac{1}{2}}
,$$
and substituting here 
(\ref{q1}) we obtain
\begin{equation}
 \omega_k \le \frac{2\cdot 12^{\frac{1}{8}}\kappa}{d_k^{\frac{9}{8}}q^{\frac{9}{8}}R^{\frac{n}{4}}}\cdot \left(\frac{2^k}{R}\right)^{\frac{3}{4}}.
\label{q2}
\end{equation}
From Lemma 9 we see that
either
$$
q\ge R^{\frac{2}{3} (n-1)}
$$
or
$$
\omega_k \ge \frac{\delta}{2q^{\frac{3}{2}}}.
$$
 From the last inequality and (\ref{q2}) we see that
$$
q^{\frac{3}{8}} = q^{\frac{3}{2}-\frac{9}{8}}\ge
 \frac{1}{4\cdot 12^{\frac{1}{8}}}\cdot\frac{\delta}{\kappa}
 \cdot  d_k^{\frac{9}{8}} R^{\frac{n}{4}}\cdot \left(\frac{R}{2^k}\right)^{\frac{3}{4}}.
$$
So in any case
$$
q^{\frac{3}{8}} \ge
\min \left( R^{\frac{1}{4} (n-1)}, \frac{1}{4\cdot 12^{\frac{1}{8}}}\cdot\frac{\delta}{\kappa}
 \cdot  d_k^{\frac{9}{8}} R^{\frac{n}{4}}\cdot \left(\frac{R}{2^k}\right)^{\frac{3}{4}}\right)=
\frac{1}{4\cdot 12^{\frac{1}{8}}}\cdot\frac{\delta}{\kappa}
 \cdot  d_k^{\frac{9}{8}} R^{\frac{n}{4}}\cdot \left(\frac{R}{2^k}\right)^{\frac{3}{4}}
$$
(to see that the minimum attains on the second element we take into account that the choice of parameters
(\ref{kkk},\ref{de}) shows that the first element in the minimum is greater than the second by the factor $R^{\frac{2}{55}}$).
Substituting the last inequality into (\ref{q2}) we obtain
$$
\omega_k \le \frac{2^7 \sqrt{12} \kappa}{d_k^{\frac{9}{2}} R^n}\cdot \left( \frac{\kappa}{\delta}\cdot \frac{2^k}{R}\right)^3.
$$

Now we must note that the number of segments of the form (\ref{form}) which intersect with   intervals
$\overline{\Delta} (A,B,C)$ corresponding to the lines from the collection  $\hbox{\got A}$
(recall that all the lines from the collection  $\hbox{\got A}$ intersect the segment $\Omega$
of the length $2\omega_k$) is 
$$
\le 
\frac{2\omega_k}{\kappa /R^{n+1}}+ 2K_k+2 \le  2^{11}R^{\frac{52}{55}}
$$
by (\ref{de},{\ref{kd}).

As for the number of   segments of the form (\ref{form}) which intersect with intervals
$\overline{\Delta} (A,B,C)$ corresponding to the lines from the collection  $\hbox{\got B}$
 we can say (Lemma 6) that  this number is
$$
\le
d_k (2K_k+2) \le 4 R^{\frac{52}{55}}
$$
by ({\ref{kd}).

In the case 2
it may happen also that a line $L(A,B,C)$ does not intersect the segment $J_n$ but the corresponding interval
$\overline{\Delta}(A,B,C)$  does intersect. But obviously such intervals can totally intersect not more than
$2K_k+2$ segments of the form  (\ref{form}).

So the total number of   segments of the form (\ref{form}) which intersect with some intervals
$\overline{\Delta} (A,B,C)$ under consideration is
$$\le
\frac{2\omega_k}{\kappa /R^{n+1}}+ (d_k+2)(2K_k+2) \le 2^{12}R^{\frac{52}{55}}.
$$
 Fundamental Lemma 1 follows as $k$ takes its values  in the interval 
$0\le k \le \log R /\log 2 $ (see
(\ref{kkk})).$\Box$

{\bf 9. Lines with large  coefficient $|A|/B$: parameter $l$.}

Here we take
an integer $l$
such that
$$
1\le l\le \frac{n}{3\lambda}
$$
 and suppose that
\begin{equation}
R^{\frac{n}{3}-\lambda (l+1)} \le B\le
R^{\frac{n}{3}-\lambda l} 
\label{elll}
\end{equation}

In this section we consider a single segment $J_{n-l} = J_{n-l}^\nu$ from the collection (\ref{coll})
with fixed  lower index $n-l$. 
 
 Let 
$$
J_n^\nu,\,\,\,\,\,1\le \nu \le T
$$ be {\it all} the segments such that
$$  
J_n^\nu \cap \overline{\Delta }(A,B,C) =\varnothing
$$
for all triples $A,B,C$  such that $ H(A,B) < R^{n-1}$ and
$$ 
J_n^\nu \subset J_{n-l}.
$$ 
Each of the segments $ 
J_n^\nu $ we divide into $R$ smaller segments
\begin{equation}\label{segmentsell}
I_{n+1}^{\nu ,\mu} ,\,\,\,1\le \nu \le T,\,\,\,
1\le \mu \le R
\end{equation}
 of equal length
$$
|I_{n+1}^{\nu ,\mu} |=
\frac{|J_n^\nu|}{R} = \frac{\kappa}{R^{n+1}}.
$$
such that
$$
J_n^\nu = \bigcup_{1\le \mu\le R} I_{n+1}^{\nu ,\mu},\,\,\,\,
1\le \nu \le T.
$$
The purpose of
  the current section is to prove that 
{\it   the number of   segments of the form
(\ref{segmentsell}) satisfying
 \begin{equation}\label{uuuu1}
I_{n+1}^{\nu,\mu}\cap \overline{\Delta}(A,B,C) \not= \varnothing
\end{equation}
for some interval  $\overline{\Delta} (A,B,C)$ with coefficients $A,B$ satisfying the  conditions (\ref{step1}) and satisfying the {\it additional} condition (\ref{elll})
is  $$\le \gamma_1 R^{\frac{52}{55}} .$$}

An admissible value for $\gamma_1$ is $\gamma_1 = 8$.

Under the conditions (\ref{step1},\ref{elll})
one has
 \begin{equation}\label{A1}
|A|\le 
 R^{\frac{n}{3}+\frac{\lambda ( l+1)}{2}}
\end{equation}
 and
 \begin{equation}\label{A1frac}
\frac{|A|}{B}\le 
 R^{\frac{3}{2}\lambda (l+1)}.
\end{equation}

In the rest part of this section we  modify lemmas 1 -  4 and 9  is the case of the inequalities (\ref{elll}).
Proofs of all lemmas below  are quite similar to the proofs of lemmas behind.

{\bf 9.1. Modified lemmata about lines intersecting a segment.}

{\bf Lemma 1$^*$.}\,\, {\it Consider a segment  $J_{n-l} = J_{n-l}^\nu$. Suppose that there exist two lines
$$
L_1= L(A_1,B_1,C_1),\,\,\, L_2= L(A_2,B_2,C_2)
$$
such that
$$
L_i \cap J_{n-l} \neq \varnothing ,\,\,\,\, i = 1,2
$$
and
$$
H(A_1,B_1), H(A_2,B_2) < R^n.
$$
 Then lines $L_1$ and $L_2$ are not parallel.}

Proof.  \,\, Lines $L_i, i = 1,2$  intersect the segment $J_{n-l} \subset \Theta$ in 
 points 
$$
\left( \theta, 
\frac{A_1\theta+C_1}{B_1}\right),\,\,\,\
\left(\theta,
\frac{A_2\theta+C_2}{B_2}\right)
$$
with $y$-coordinates
$$
\frac{A_1\theta+C_1}{B_1},\,\,\,\
\frac{A_2\theta+C_2}{B_2},\,\,\,\,
\left|\frac{A_1\theta+C_1}{B_1}-
\frac{A_2\theta+C_2}{B_2}\right|\le |J_{n-l}|.
$$
Suppose these lines to be parallel.
Then
$$
\frac{A_1}{B_1} =\frac{A_2}{B_2}
$$
and by making use of (\ref{elll}) we have
$$
\frac{\kappa}{R^{n-l}} <
\frac{R^{\frac{n}{3}+(2\lambda -1) l}}{R^{n-l}}=
\frac{R^{2\lambda l}}{R^{\frac{2n}{3}}}\le
\frac{1}{B_1B_2} \le 
\left|
\frac{C_1}{B_1} -\frac{C_2}{B_2}
\right| =
\left|
\frac{A_1\theta+C_1}{B_1} -
\frac{A_2\theta+C_2}{B_2}\right| \le |J_{n-l}| = \frac{\kappa}{R^{n-l}}
$$
(here we use the inequality $\kappa<1< R^{\frac{n}{3}+(2\lambda -1) l}$) and this is a contradiction. $\Box$

We do not need any changes in Lemma 2.
 But in the case $l\ge 1$ simple application of Lemma 1 gives a strong inequality.
This inequality 
we formulate as

{\bf Lemma 2$^*$}\,\,\,{\it
Suppose that two lines 
$$
L_1= L(A_1,B_1,C_1),\,\,\, L_2= L(A_2,B_2,C_2),\,\,\,\,
L_i \cap J_{n-l} \neq \varnothing ,\,\,\,\, i = 1,2
$$
satisfy 
$$
H(A_i,B_i) \le R^n,\,\,\,\, i = 1,2
.$$
Suppose the additional condition (\ref{elll}) to be valid.
Then
$$
|q\theta - p|\le \kappa R^{-\frac{n}{3} - (2\lambda - 1)l}.
$$}

Proof.\,\,\,  We should take in Lemma 2 
$I=J_{n-l}$ 
and combine the conclusion (i) with  (\ref{elll}).$\Box$

{\bf Lemma 3$^*$.}\,\, {\it
All the lines $ L = L(A,B,C)$ 
such that 
$$
L(A,B,C) \cap J_{n-l} \neq \varnothing,
$$
$$
H(A,B) < R^n
$$
satisfying the addditional condition (\ref{elll})
have a single common point.}

Proof.  

The proof is quite close to the proof of Lemma 3. 
  From Lemma 1$^*$  it follows that any two lines  intersecting $J_{n-l}$  have a common point.
Suppose that we have three lines 
$$
L_i=
L(A_i,B_i,C_i),\,\,\,\, i = 1,2,3
$$
intersecting $J_{n-l}$  which satisfy the conditions of Lemma 3$^*$ but do not have  a common point.
Then
  (by taking  $(\theta,\ \xi)$ to be the middle of $J_{n-l}$ and 
   $(\theta, Y_i) = J_{n-l}\cap L_i$ we see that
$$
|A_i\theta -B_i\xi +C_i|=
B_i |\xi - Y_i| \le  B_i \cdot \frac{\kappa}{2R^{n-l}}
$$
for every $i =1,2,3$) we have
$$1\le
|D| =|\,
\left|
\begin{array}{ccc}
A_1&B_1&C_1\cr
A_2&B_2&C_2\cr
A_3&B_3&C_3
\end{array}
\right|\,| \le
 3\cdot 2AB \cdot \frac{\kappa B}{2R^{n-l}} = \frac{3\kappa AB^2}{R^{n-l}},
$$
where 
$$A = \max_{i=1,2,3}  |A_i| < R^{\frac{n}{3}  +\frac{(l+1)\lambda}{2}}
,\,\,\, B = \max_{i=1,2,3} B_i
\le R^{\frac{n}{3}-\lambda l}
$$
(the first  inequality here follows from inequalities (\ref{elll}) as
$ |A_i|^2 \le \frac{R^n}{B_i} \le\frac{R^n}{R^{\frac{n}{3}-(l+1)\lambda }} = R^{\frac{2}{3}n + (l+1)\lambda}$).

   Recall that we suppose the condition  (\ref{00k}) to be valid and $\lambda$ satisfies (\ref{lala}). So
$$
1\le |D| < 3\kappa R^{\frac{(l+1)\lambda}{2}- (2\lambda -1)l} < 1.$$ 
This is not possible and lemma is proved.$\Box$

Now we  suppose that all the lines intersecting the segment $J_{n-l}$ and satisfying $ H(A,B) < R^{n-1}$ 
and the additional condition  (\ref{elll})  pass through a single point
$$
P=\left(\frac{p}{q},\frac{r}{q}\right).
$$
Put
\begin{equation}\label{sigmanew}
\sigma (l) = \frac{\kappa R^l}{|q\theta -p|}.
\end{equation}

{\bf Lemma 4$^*$.}\,\,\,{\it  Consider two lines $L_i = L(A_i,B_i,C_i), i = 1,2$.
Suppose that $A_1,B_1,A_2,B_2$ satisfy 
$$
H(A_1,B_1), H(A_2,B_2) <R^n.
$$
 Suppose that 
both lines $L_1,L_2$ intersect the segment $J_{n-l}$. Suppose that 
$$
B_1\ge B_2.
$$
 
Then with $\sigma (l)$ defined in (\ref{sigmanew}) one has
\begin{equation}\label{ppo1new}
B_1\le {\sigma}(l),\,\,\, A_1^2\le {\sigma}(l) \,  B_1.
\end{equation}}

Proof.\,\,\,
By Lemma 2  (statement (i)) we  have
$$
|q\theta - p| \le \frac{\kappa}{R^{n-l}}\cdot B_1B_2
\le
\frac{\kappa}{R^{n-l}}\cdot B_1^2.
$$
As in Lemma 4 we see that
$$
\max\left(
\frac{A_1^2}{B_1}, B_1
\right)=
\frac{H(A_1,B_1)}{B_1^2} \le
 \frac{\kappa}{|q\theta -p|}\cdot {R^l} = {\sigma (l)}
$$
 and Lemma 4$^*$ follows.$\Box$

Put
\begin{equation}\label{veel}
V(l) := 
\left(\frac{(\sigma (l))^3}{R^{n-1}}\right)^{\frac{1}{4}}.
\end{equation}

{\bf Corollary 1.}\,\,\,{\it 
Suppose that the conditions of 
 Lemma 4$^*$ are satisfied and in addition we have (\ref{step1}). Then
\begin{equation}\label{corolla}
\frac{|A_1|}{B_1} \le
 V(l).
\end{equation}}

Proof.\,\, Apply Lemma 5  with
$\sigma = \sigma (l), W = R^{n-1}$ .$\Box$

{\bf Corollary 2.}\,\,\,{\it 
 Let
\begin{equation}\label{ul}
L_1,L_2,..., L_M , \,\,\,\,  L_j = L(A_j,B_j,C_j)
\end{equation}
be all   the lines intersecting $J_{n-l}$ and satisfying (\ref{step1}).
Then  for all j from the interval $ 1\le j \le M$  but one possible exception one has
\begin{equation}\label{corolla2}
\frac{|A_j|}{B_j} \le
V(l)  .
\end{equation}} 
 
Proof.  \,\, Among the collection (\ref{ul}) we have a line with the minimal coefficient $B_j$. By 
(\ref{corolla}) of Corollary 1 we see that all other lines satisfy (\ref{corolla2}).
$\Box$

{\bf 9.2. Collections  $\hbox{\got A}_l$   and $\hbox{\got B}_l $.}

In the sequel we  suppose that $M\ge 2$.
We divide the collection of lines (\ref{ul}) into two subcollections.
Collection $\hbox{\got B}_l $ consists of   only one line   with the minimal value of $B$.
So
\begin{equation}\label{odin}
\#\hbox{\got B}_l =1.
\end{equation}
All other lines form the collection $\hbox{\got A}_l $.
By the arguments form the proof of Corolary 2 we see that for any $L_j$ from the collection $\hbox{\got A}_l $
we have (\ref{corolla2}).
So all these lines intersect the segment $\Theta$ in the points of the segment
$$\Omega (l) =\left[\frac{r}{q}-\omega (l),
\frac{r}{q}+\omega (l)\right],
$$
where
\begin{equation}\label{amegaell}
\omega (l) = \left|\theta -\frac{p}{q}\right|\cdot V(l) =
\frac{|q\theta-p| (\sigma (l))^{\frac{3}{4}} }{qR^{\frac{n-1}{4}}}=
\frac{\kappa^{\frac{3}{4}}|q\theta-p|^{\frac{1}{4}}  R^{\frac{3l+1}{4}}}{qR^{\frac{n}{4}}}
\end{equation}
by the definitions of $\sigma (l)$ and $V(l)$ (see (\ref{sigmanew},\ref{veel})).
We apply Lemma 2$^*$
to deduce from (\ref{amegaell}) the inequality
\begin{equation}\label{amegaell1}
\omega (l) 
\le \frac{\kappa}{qR^{\frac{n}{3}}} \cdot
R^{-\frac{\lambda l}{2} + \frac{4l+1}{4}}.
\end{equation}

{\bf 9.3. Collection  $\hbox{\got A}_l$: lower bound for $q$ and its application.}

We deal with the situation $l\ge 1$. In this case the consideration of the collection $\hbox{\got A}_l$ is much more simple.
The only thing  what we need is an analog of Lemma 9 and its corollary for the lower bound of $q$.

{\bf Lemma 9$^*$.}\,\,\,{\it
Suppose that 
 \begin{equation}\label{qu1}
q< R^{\frac{2}{3} (n-l-1)}.
\end{equation}
Then
 for the value $\omega (l)$ defined in (\ref{amegaell}) one has
$$
\omega (l) \ge \frac{\delta}{2q^{\frac{3}{2}}}.
$$}

Proof.\,\,\, Similarly to the proof of Lemma 9 we find integers $A,B,C$ such that $ (A,B)\neq (0,0)$  and
$$
Ap-Br+Cq=0
,\,\,\,
\max (|A|,B) \le q^{\frac{1}{2}},\,\,\, B\ge 0.
$$
Then  
$$
\left| A\theta -B\frac{r}{q} +C\right|=
\left| A\left(\theta-\frac{p}{q}\right) +\frac{Ap-Br+Cq}{q}\right| =
\left| A\left(\theta-\frac{p}{q}\right) \right|\le \frac{|q\theta- p|}{q^{\frac{1}{2}}},
$$
From the condition (\ref{qu1}) we have
$$
R^{\frac{n}{3}} \ge q^{\frac{1}{2}} R^{\frac{l+1}{3}}.
$$
So we take into account Lemma 
2$^*$   to see that
$$
\left| A\theta -B\frac{r}{q} +C\right|
\le
\frac{|q\theta- p|}{q^{\frac{1}{2}}}
\le
\frac{\kappa R^{-\frac{n}{3}- (2\lambda -1)}}{q^{\frac{1}{2}}} \le
\frac{\kappa}{   q} \cdot R^{-(2\lambda -1) - \frac{l+1}{3}}.
$$
As
$$
\frac{\delta}{q}\le \frac{\delta}{A^2}\le
|A\theta +C| 
$$
by (\ref{0})
and
 $\delta >  {\kappa}\cdot R^{-(2\lambda -1) - \frac{l+1}{3}}$ we have $B>0$.
As
$$
\max (|A|,B) < R^{\frac{n-l-1}{3}}
$$
it follows that
$$
\overline{\Delta}  (A,B,C) \cap J_{n-l} =\varnothing.
$$
By following all the arguments of the proof of Lemma 9 we see that
$$
\omega (l) + \frac{\kappa}{B  q}\cdot R^{-(2\lambda -1) - \frac{l+1}{3}}
\ge \frac{\delta}{Bq}.
$$
 As  $\lambda >3$ and $\kappa ={\delta}R^{\frac{6}{5}}$   we have
$$\frac{\delta}{2} > {\kappa}\cdot R^{-(2\lambda -1) - \frac{l+1}{3}}.$$
Lemma 9 $^*$ follows.$\Box$

{\bf Corollary 1.}\,\,\,
{\it The following inequality is valid:
$$
q\ge R^{\frac{2}{3} (n-l-1)}.$$
}

Proof.\,\,
Suppose that (\ref{qu1}) is valid. Then
by Lemma 9$^*$ we have
$$
\omega (l) \ge \frac{\delta}{2q^{\frac{3}{2}}}.
$$
Combining this inequality with (\ref{amegaell1}) we have
$$
\frac{\delta}{2q^{\frac{3}{2}}}\le
\frac{\kappa}{qR^{\frac{n}{3}}} \cdot
R^{-\frac{\lambda l}{2} + \frac{4l+1}{4}}.
$$
Hence
$$
q \ge \frac{1}{4}\, \left(\frac{\delta}{\kappa}\right)^2 R^{\frac{2}{3}n +\lambda l -\frac{4l+1}{2}}=
\frac{1}{4}\cdot R^{\frac{2}{3}n + (\lambda - 2)l  - \frac{29}{10}}> R^{\frac{2}{3} (n-l-1)},
$$
as $\lambda > 4$ and
$$
(\lambda - 2)l  - \frac{29}{10} > 2l - \frac{29}{10}\ge - \frac{9}{10}> -\frac{4}{3}\ge \frac{2}{3} (-l-1).$$
Corollary 1 is proved.$\Box$

{\bf Corollary 2.}\,\,\,
{\it In the case $l \ge 1$ we have  the following upper bound:
$$ 
\omega (l) \le \frac{\kappa}{R^n}\cdot 
{R^{-\frac{\lambda l}{2} +\frac{20l+11}{12}}}.
$$}

Proof.\,\,\, Apply   (\ref{amegaell1}) and  the inequality of  Corollary 1.$\Box$

{\bf 9.4. The second fundamental lemma.}

Here we prove the following

{\bf Fundamental Lemma 2.}\,\,\,{\it
Let $l \ge 1$.
Suppose we have a segment $J_{n-l} $.
Then the number of segments $I_{n+1}^{\nu,\mu}$ of the form (\ref{segmentsell}) which   intersect with some interval
$$
\overline{\Delta}(A,B.C)
$$
with $A,B$ satisfying (\ref{step1}) and (\ref{elll}) is
$$
\le 8 R^{\frac{52}{55}}.
$$
}

Proof.

First  of all we suppose that $M\ge 2$ (otherwise  there
exists
only one line $L_1$ under consideration and we may  use the san=me arguments as
for the collection $\hbox{\got B}_l$, see below).

1. Lines from the collection $\hbox{\got A}_l$ intersect the segment $\Theta$. The points of intersection
belong to the  segment  $\Omega_l$
of the length $2\omega (l)$ satisfying upper bound  (\ref{amegaell1}).
For $L(A,B,C) \in  \hbox{\got A}_l$ one has
$$
|\overline{\Delta} (A,B,C) | =\frac{2\delta}{H(A,B)}\le\frac{2\delta}{R^{n-1}}.
$$
So the  number of segments 
$
I_{n+1}^{\nu,\mu}$ of the form (\ref{segmentsell}) which  intersect
with intervals $\overline{\Delta} (A,B,C)$ corresponding to the collection  $\hbox{\got A}_l$
is less or  equal than
$$
\frac{2\omega (l) +\max |\overline{\Delta} (A,B,C) |}{\kappa/R^{n+1}}+ 2 \le
2R^2\cdot \frac{\delta}{\kappa} + 2 R^{-\frac{\lambda l}{2}+ \frac{20l+23}{12}}+2\le 
2R^2\cdot \frac{\delta}{\kappa} + 2 R^{-\frac{\lambda }{2}+ \frac{43}{12}}+2\le 4 R^{\frac{52}{55}}+2
$$
as $\frac{\delta}{\kappa} = R^{-\frac{6}{5}}$ and $\lambda  = \frac{1741}{330}$.

2. The 
  number of segments 
$
I_{n+1}^{\nu,\mu}$ of the form (\ref{segmentsell}) which  intersect
with intervals $\overline{\Delta} (A,B,C)$ corresponding to the collection  $\hbox{\got B}_l$
is less or  equal than
$$
\frac{ 2\delta /R^{n-1}}{\kappa/R^{n+1}}+ 2 \le  2 R^{\frac{4}{5}}+2.
$$

Also we must take into account that a line $L(A,B,C)$ may not intersect the segment $J_{n-l}$ but the corresponding interval
$\overline{\Delta}(A,B,C)$  may intersect it. But obviously such intervals can totally intersect not more than
$2R^{\frac{4}{5}}+2$ segments of the form  (\ref{segmentsell}).

The second Fundamental Lemma follows.$\Box$

{\bf 10.  Proof  of Proposirion 1.}

We apply Fundamental Lemmas 1 and 2.  Arguments below are close to those from Peres-Schlag's method (see \cite{PS}).

Recall that we denote by $T_n$ the total number of   segments $J_n^\nu$.

By Fundamental Lemmas 1,2 we see that
$$
T_{n+1} \ge T_n\cdot R    - T_n\cdot  2^{13}R^{\frac{52}{55}}\log R
-\sum_{l=1}^{[n/3\lambda ]} T_{n-l} \cdot 8R^{\frac{52}{55}}
$$
or
$$
T_{n+1} \ge T_n
\left(
 R    -   2^{13}R^{\frac{52}{55}}\log R
-\sum_{l=1}^{[n/3\lambda ]} \frac{T_{n-l}}{T_n} \cdot 8R^{\frac{52}{55}}\right).
$$
We see by induction that
$$
T_{n+1} \ge T_n \cdot (R-  2^{14}R^{\frac{52}{55}}\log R)
$$
or
$$
T_n \ge (R-  2^{14}R^{\frac{52}{55}}\log R)^{n-1}.
$$
In fact as $R\ge 2^{422}$ this inequality proves that $T_n >0$ for every $n$.
It means that
$$
\bigcap_{n \in \mathbb{N} }\,\,\bigcup_{1\le \nu_n\le T_n}\,\, J_n^{\nu_n} \neq \varnothing .
$$
By putting
$R =2^{422} $ we prove Proposition 1. $\Box$
 
\vskip+1.0cm

\vskip+1.0cm

author: Nikolay G. Moshchevitin;

e-mail: moshchevitin@rambler.ru 

\begin{thebibliography}{99}
\bibitem{SCH}  W. M. Schmidt,\,\,\,
 Open problems in Diophantine approximations. //  "Approximations Diophantiennes et nombres transcendants" Luminy, 1982, Progress in
 Mathematics,
 Birkh\"auser (1983), p.271 - 289.

\bibitem{PSV} D.Badziahin, 
A.Pollington, S.Velani,\,\,\,
On a problem in simultaneous Diophantine approximations: Schmidt's conjecture. //
Preprint, available at arXiv:1001.2694v2 10Mar2010, 


\bibitem{PS} Y. Peres, W. Schlag, \,\,\, Two Erd\"os problems on lacunary sequences: chromatic numbers
    and Diophantine approximations.  // Preprint, available at:
 arXiv:0706.0223v1  1Jun2007.
  \end{thebibliography}
\end{document}